\documentclass[reqno,centertags,12pt]{amsart}

\usepackage{amssymb,amsmath,amsfonts,amssymb} 
-\marginparwidth \oddsidemargin -9mm \evensidemargin \oddsidemargin
\usepackage{latexsym}
\advance\hoffset by 5mm

\def\@abssec#1{\vspace{.05in}\footnotesize \parindent .2in
{\bf #1. }\ignorespaces}
\newtheorem{theorem}{Theorem}[section]
\newtheorem{thm}[theorem]{Theorem}
\newtheorem{lemma}[theorem]{Lemma}

\newtheorem{definition}[theorem]{Definition}

\def \Rm {\mathbb R}
\def \Tm {\mathbb T}

\def \Zm {\mathbb Z}
\def\cU{\mathcal U}
\def\cS{\mathcal S}

\allowdisplaybreaks \numberwithin{equation}{section}

\title{A variation on a theme of Caffarelli and Vasseur}
\author{Alexander Kiselev and Fedor Nazarov} 
\dedicatory{Dedicated to Nina Nikolaevna Uraltseva}
\thanks{Department of
Mathematics, University of Wisconsin, Madison, WI 53706, USA;
email: kiselev@math.wisc.edu, nazarov@math.wisc.edu}

\begin{document}


\begin{abstract}
Recently, using DiGiorgi-type techniques, Caffarelli and Vasseur
\cite{CV} showed that a certain class of weak solutions to the drift
diffusion equation with initial data in $L^2$ gain H\"older
continuity provided that the BMO norm of the drift velocity is
bounded uniformly in time. We show a related result: a uniform bound
on BMO norm of a smooth velocity implies uniform bound on the
$C^\beta$ norm of the solution for some $\beta >0.$ We use
elementary tools involving control of H\"older norms using test
functions. In particular, our approach offers a third proof of the
global regularity for the critical surface quasi-geostrophic (SQG)
equation in addition to \cite{KNV} and \cite{CV}.
\end{abstract}

\maketitle

\section{Introduction}\label{intro}

In the preprint \cite{CV}, Caffarelli and Vasseur proved that
certain weak solutions of the drift diffusion equation with
$(-\Delta)^{1/2}$
dissipation 
gain H\"older regularity provided that the velocity $u$ is uniformly
bounded in the BMO norm. The proof uses DiGiorgi-type iterative
techniques. The goal of this paper is twofold. First, we wanted to
provide additional intuition for the Caffarelli-Vasseur theorem by
presenting an elementary proof of a related result. Secondly, we
think that, perhaps, the method of this paper may prove useful in
other situations.

Everywhere in this manuscript, our setting for the space variable
will be $d-$dimensional torus, $\Tm^d.$ Equivalently, we may think
of the problem set in $\Rm^d$ with periodic initial data. With the
latter interpretation in mind, let us recall the definition of the
BMO norm:
\begin{equation}\label{BMO}
\|f\|_{BMO} = {\rm sup}_{B \in \Rm^d} \frac{1}{|B|} \int\limits_B
|f(x)-\overline{f}_B|\,dx.
\end{equation}
Here $B$ stands for a ball in $\Rm^d,$ $|B|$ for its volume and
$\overline{f}_B$ for the mean of the function $f$ over $B.$
\begin{thm}\label{main1}
Assume that $\theta(x,t),$ $u(x,t)$ are $C^\infty(\Tm^d \times
[0,T])$ and such that
\begin{equation}\label{ps}
\theta_t = (u \cdot \nabla) \theta - (-\Delta)^{1/2} \theta
\end{equation}
holds for any $t \geq 0.$ Assume that the velocity $u$ is divergence
free and satisfies a uniform bound $\|u(\cdot,t)\|_{\rm{BMO}} \leq
B$ for $t \in [0,T].$ Then there exists $\beta = \beta(B,d)>0$ such
that
\begin{equation}\label{kr}
\|\theta(x,t)\|_{C^\beta(\Tm^d)} \leq C(B,\theta(x,0))
\end{equation}
for any $t \in [0,T].$
\end{thm}
\it Remark. \rm In fact, we get control of H\"older continuity in
terms of just $L^1$ norm of $\theta_0$ if we are willing to allow
time dependence in \eqref{kr}. Namely, the following bound is also
true:
\begin{equation}\label{kr12}
\|\theta(x,t)\|_{C^\beta(\Tm^d)} \leq C(B,\|\theta(x,0)\|_{L^1}){\rm
min}(1,t)^{-d-\beta}.
\end{equation}

Thus uniform bound on the BMO norm of $u$ implies uniform bound on
a certain H\"older norm of $\theta.$ The dimension $d$ is
arbitrary.

Our result is different from \cite{CV}. 
For one thing, \cite{CV} contains a local regularization version,
which we do not attempt here. However our proof is simpler, and is
quite elementary. At the expense of extra technicalities, it can
be extended to more general settings.

Theorem~\ref{main1} can be used to give a third proof of the global
regularity of the critical surface quasi-geostrophic (SQG) equation,
which has been recently established in \cite{KNV} and \cite{CV}. We
discuss this in Section~\ref{SQG}. Throughout the paper, we will
denote by $C$ different constants depending on the dimension $d$
only.

\section{Preliminaries}\label{prem}


First, we need an elementary tool to characterize
H\"older-continuous functions. 
Define a function $\Omega(x)$ on $\Tm^d$ by
\begin{equation}\label{confun}
\Omega(x) = \left\{ \begin{array}{ll} |x|^{1/2}, & |x|<1/2 \\
\frac{1}{\sqrt{2}}, & |x| \geq 1/2 \end{array} \right.
\end{equation}
(thinking of $\Rm^d$ picture, $\Omega$ is defined as above on a unit
cell and continued by periodicity). Let $A>1$ be a parameter to be
fixed later.
\begin{definition}
We say that a $C^\infty$ function $\varphi$ defined on $\Tm^d$
belongs to $\cU_r(\Tm^d)$ if
\begin{eqnarray}\label{philinf}
\|\varphi(x)\|_{L^\infty} \leq \frac{A}{r^d} & \\
\label{phimz} \int_{\Tm^d} \varphi(x)\,dx =0 &  \\
\label{philone} \|\varphi(x)\|_{L^1} \leq 1 & \\
\label{phiconc} \int_{\Tm^d} |\varphi(x)| \Omega(x-x_0)\,dx \leq
r^{1/2} & \,\,\,\,{\rm for}\,\,\,{\rm some}\,\,\,x_0 \in \Tm^d.
\end{eqnarray}
\end{definition}
Observe that the classes $\cU_r$ are invariant under shift. We will
write $f(x) \in a\cU_r(\Tm^d)$ if $f(x)/a \in \cU_r(\Tm^d).$
The choice of the exponent $1/2$ in \eqref{confun} and
\eqref{phiconc} is arbitrary and can be replaced with any positive
number less than $1$ with the appropriate adjustment of the range of
$\beta$ in Lemma~\ref{Uhold} below.

The classes $\cU_r$ can be used to characterize H\"older spaces as
follows. Let us denote
\begin{equation}\label{Hnorm}
\|f\|_{C^\beta(\Tm^d)} = {\rm sup}_{x,y \in
\Tm^d}\frac{|f(x)-f(y)|}{|x-y|^\beta},
\end{equation}
omitting the commonly included on the right hand side
$\|f\|_{L^\infty}$ term. The seminorm \eqref{Hnorm} is sufficient
for our purposes since $\theta$ remains bounded automatically, and
moreover we could without loss of generality restrict consideration
to mean zero $\theta,$ invariant under evolution, for which
\eqref{Hnorm} is equivalent to the usual H\"older norm.
\begin{lemma}\label{Uhold}
A bounded function $\theta(x)$ is in $C^\beta(\Tm^d),$ $0 < \beta <
1/2,$ if and only if there exists a constant $C$ such that for every
$0 < r \leq 1,$
\begin{equation}\label{holdcon}
\left| \int_{\Tm^d} \theta(x)\varphi(x)\,dx \right| \leq C r^\beta
\end{equation}
for all $\varphi \in \cU_r.$ Moreover,
\begin{equation}\label{holdcon121}
\|\theta\|_{C^\beta(\Tm^d)}\leq C(\beta) {\rm sup}_{\varphi \in
\cU_r,\,\,0 <r \leq 1} r^{-\beta}\left|\int_{\Tm^d}
\theta(x)\varphi(x)\,dx \right|.
\end{equation}
\end{lemma}
\noindent \it Remark. \rm The lemma holds for each fixed $A$ in
\eqref{philinf}. It will be clear from the proof that the constant
$C$ in \eqref{holdcon121} does not depend on $A$ provided that $A$
was chosen sufficiently large.
\begin{proof}
Assume first that $\theta \in C^\beta.$ Consider any $\varphi \in
\cU_r,$ and 
observe that
\[ \int_{\Tm^d} \theta(x)\varphi(x)\,dx = \int_{\Tm^d}
(\theta(x)-\theta(x_0))\varphi(x)\,dx \leq C\int_{\Tm^d}
|x-x_0|^\beta |\varphi(x)|\,dx. \]
Using H\"older inequality, we get \begin{eqnarray}\nonumber
\int_{\Tm^d} |x-x_0|^\beta |\varphi(x)|\,dx \leq \left(\int_{\Tm^d}
|\varphi(x)|\,dx\right)^{1-2\beta}\left(\int_{\Tm^d} |x-x_0|^{1/2}
|\varphi(x)|\,dx\right)^{2\beta} \\ \label{th1} \leq
C\left(\int_{\Tm^d}
|\varphi(x)|\,dx\right)^{1-2\beta}\left(\int_{\Tm^d} \Omega(x-x_0)
|\varphi(x)|\,dx\right)^{2\beta}
\end{eqnarray}
Due to \eqref{philone} and \eqref{phiconc}, the right hand side of
\eqref{th1} does not exceed $r^\beta.$

For the converse, consider a periodization $\theta_p$ of $\theta$ in
$\Rm^d.$ Recall a well known characterization of H\"older continuous
functions in $\Rm^d$ (see e.g. \cite{Stein}):
\begin{equation}\label{holdcrit}
 \theta \in C^\beta \,\,\Leftrightarrow\,\,\|\theta\|_{L^\infty}
 \leq Q,\,\,\,\|\Delta_j(\theta)\|_\infty \leq Q2^{-\beta
j},\,\,\,\forall j.
\end{equation}
Moreover, if the right hand side of \eqref{holdcrit} is satisfied,
then $\|\theta\|_{C^\beta} \leq CQ.$

Here $\Delta_j$ are the Littlewood-Paley projections:
\[ \Delta_j(\theta) = \theta \ast \Psi_{2^{-j}}, \]
where $\Psi_t(x) = t^{-d}\Psi(x/t),$ and $\widehat{\Psi}(\xi) =
\eta(\xi)-\eta(2\xi),$ with $\eta \in C_0^\infty,$ $0 \leq \eta(\xi)
\leq 1,$ $\eta(\xi)=1$ if $|\xi|\leq 1$ and $\eta(\xi)=0,$ $|\xi|
\geq 2.$ Observe that $\Psi$ is in the Schwartz class $\cS,$
$\int_{\Rm^d} \Psi \,dx =0,$ $\int_{\Rm^d} |\Psi(x)|\,dx \leq C,$
$\int_{\Rm^d} |x|^{1/2}|\Psi(x)|\,dx \leq C$ and
$\|\Psi\|_{L^\infty} \leq C.$ Let us define
\[ \widetilde{\Psi}_j(x) =
c\sum\limits_{n \in \Zm^d} \Psi_{2^{-j}}(x+n), \] then
$\widetilde{\Psi}_j(x) \in \cU_{2^{-j}}(\Tm^d)$ if $c$ is
sufficiently small (independently of $j$). Moreover,
\begin{equation}\label{aux11}
\int_{\Tm^d} \theta(x) \widetilde{\Psi}_j(x-y)\,dx =
c\int_{\Rm^d} \theta_p(x) \Psi_{2^{-j}}(x-y)\,dx.
\end{equation}
By assumption, the left hand side in \eqref{aux11} does not exceed
$Q2^{-j\beta}.$ Thus, by the criterion \eqref{holdcrit}, $\theta_p$
is $C^\beta$ and so is $\theta.$ The remark after \eqref{holdcrit}
implies that \eqref{holdcon121} is true.
\end{proof}

The proof of Theorem~\ref{main1}, which we will outline in the
beginning of the next section, relies on transfer of evolution on
the test function. Here is an elementary lemma that allows us to
do that. Let $\varphi^t(x,s)$ be the solution of
\begin{equation}\label{phieq}
\varphi^t_s = -(u(x,t-s)\cdot\nabla)\varphi^t -
(-\Delta)^{1/2}\varphi^t,\,\,\,\,\varphi^t(x,0)=\varphi(x).
\end{equation}
\begin{lemma}\label{evoltran}
Let $\theta_0, \varphi \in C^\infty(\Tm^d),$ and let $\theta(x,t)$
be the solution of \eqref{ps} with $\theta(x,0)=\theta_0(x).$ Then
we have
\[ \int_{\Tm^d} \theta(x,t) \varphi(x)\,dx = \int_{\Tm^d} \theta_0(x)
\varphi^t(x,t)\,dx. \]
\end{lemma}
\begin{proof}
We claim that for $0 \leq s \leq t,$ the expression
\begin{equation}\label{thphiex}
 \int_{\Tm^d} \theta(x, t-s) \varphi^t(x,s)\,dx
\end{equation}
remains constant. 
A direct computation using \eqref{phieq}, \eqref{ps} and the fact
that $u$ is divergence free shows that the $s$-derivative of
\eqref{thphiex} is zero. Substituting $s=0$ and $s=t$ into
\eqref{thphiex} proves the lemma.
\end{proof}

\section{The proof of the main result}\label{bs}

Let us outline our plan for the proof of Theorem~\ref{main1}.
Conceptually, the proof is quite simple: integrate the solution
against a test function from $\cU_r,$ transfer the evolution on the
test function and prove estimates on the test function evolution. 

The key to the proof of Theorem~\ref{main1} is the following result.
\begin{thm}\label{key}
Let $v(x,s)\in C^\infty (\Tm^d \times [0,T])$ be divergence free
d-dimensional vector field, and let $\psi(x,s)$ solve
\begin{equation}\label{psieq}
\psi_s = -(v \cdot \nabla)\psi - (-\Delta)^{1/2}\psi,
\,\,\,\psi(x,0)=\psi(x).
\end{equation}
Assume that
\[ {\rm max}_{s \in [0,T]} \|v(\cdot, s)\|_{BMO} \leq B.  \]
Then the constant $A=A(B,d)$ in \eqref{philinf} can be chosen so
that the following is true.

Suppose $\psi \in \cU_r(\Tm^d),$ $0<r \leq 1.$ Then there exist
constants $\delta$ and $K>0,$ which depend only on $B$ and dimension
$d,$ such that
\begin{equation}\label{testevol11}
\psi(x,s) \in \left(\frac{r}{r+Ks}\right)^{\delta/K}
\cU_{r+Ks}(\Tm^d)
\end{equation}
if $r+Ks \leq 1$ and $\psi(x,s) \in r^{\delta/K} \cU_1(\Tm^d)$
otherwise.
\end{thm}

Let us assume Theorem~\ref{key} is true and prove
Theorem~\ref{main1}.

\begin{proof}[Proof of Theorem~\ref{main1}]
Let $\beta = \delta/K.$ Since $\theta(x,0)$ is smooth, we have
\begin{equation}\label{ini}
\left| \int_{\Tm^d} \theta(x,0) \varphi(x)\,dx \right| \leq
C(\theta(x,0))r^\beta
\end{equation}
for all $\varphi(x) \in \cU_r(\Tm^d),$ $0<r \leq 1.$ But by
Lemma~\ref{evoltran},
\[ \int_{\Tm^d} \theta(x,t)\varphi(x)\,dx = \int_{\Tm^d} \theta(x,0)
\varphi^t(x,t)\,dx. \] By Theorem~\ref{key}, $\varphi^t(x,t)$
belongs to  $\left(\frac{r}{r+Kt}\right)^{\beta} \cU_{r+Kt}(\Tm^d)$
if $r+Kt \leq 1,$ and to $r^{\beta}\cU_1(\Tm^d)$ otherwise. Then
\eqref{ini} implies that \begin{equation}\label{finb14} \left|
\int_{\Tm^d} \theta(x,t)\varphi(x)\,dx \right| \leq
C(\theta(x,0))r^\beta,
\end{equation} for all $\varphi(x) \in \cU_r(\Tm^d),$ $0<r \leq 1.$
\end{proof}
 Observe that $C(\theta(x,0))$ will depend only on the $L^1$ norm
of $\theta(x,0)$ if we are willing to allow time dependence in
\eqref{finb14}:
\begin{eqnarray*}
\left| \int_{\Tm^d} \theta(x,t)\varphi(x)\,dx \right| = \left|
\int_{\Tm^d} \theta(x,0) \varphi^t(x,t)\,dx \right| \leq
\|\theta(x,0)\|_{L^1}\|\varphi^t(x,t)\|_{L^\infty} \\ \leq
A(B,d)r^\beta {\rm min}(1,r+Kt)^{-d-\beta} \leq
C(B,d,\|\theta(x,0)\|_{L^1}){\rm min}(1,t)^{-d-\beta}r^\beta.
\end{eqnarray*}
This proves the bound \eqref{kr12} in the Remark after
Theorem~\ref{main1}.

Thus it remains to prove Theorem~\ref{key}.

\section{The Evolution of the Test Function}\label{etf}


The proof of Theorem~\ref{key} is based on the following lemma,
which looks at what happens over small time increments.

\begin{lemma}\label{keylemma}
Under assumptions of Theorem~\ref{key}, we can choose $A=A(B,d)$ so
that the following is true. There exist positive $\delta,$ $K$ and
$\gamma$ (dependent only on $B$ and $d$) such that for all $0 \leq s
\leq \gamma r,$ if $\psi(x,0) \in \cU_r(\Tm^d),$ $0<r \leq 1,$ then
\begin{equation}\label{finreslem}
\psi(x,s) \in \left(1-\frac{\delta s}{r}\right) \cU_{r+Ks}(\Tm^d).
\end{equation}
The estimate \eqref{finreslem} is valid as long as $r+Ks \leq 1;$
otherwise the solution just remains in $\cU_1.$
\end{lemma}
\begin{proof}
We have to check four conditions. First, the equation for $\psi$
preserves the mean zero property, so that $\int_{\Tm^d}
\psi(x,s)\,dx =0$ for all $s.$

Next, let us consider {\bf the  $L^\infty$ norm}. 
Set $M(s)= \|\psi(\cdot, s)\|_{L^\infty}.$ Consider any point
$x_0$ where the maximum or minimum value is achieved. Without any
loss of generality, we can assume $x_0=0,$ $\psi(0,s)=M(s).$ Then
\begin{equation}\label{maxdev}
\partial_s\psi(0,s) = -(-\Delta)^{1/2}\psi(0,s) = C \sum\limits_{n \in Z^d}
\int_{\Tm^d} \frac{\psi(y,s) - M(s)}{|y+n|^{d+1}}\,dy.
\end{equation}
Here we used the well-known  formula for the fractional Laplacian
(see e.g. \cite{CC}). 
Since $\|\psi(\cdot,s)\|_{L^1(\Tm^d)} \leq 1$ (see the argument
below on the $L^1$ norm monotonicity), it is clear that the
contribution to the right hand side of \eqref{maxdev} from the
central period cell is maximal when $\psi(y)$ is the characteristic
function of a ball of radius $cM(s)^{-1/d}$ centered at the origin.
This gives us the estimate
\begin{equation}\label{infin1}
\partial_s\psi_s(0,s) \leq -C\int_{cM(s)^{-1/d}}^{r^{-1}}
M(s)|y|^{-d-1}\,dy \leq -C_1M(s)^{\frac{d+1}{d}}+C_2rM(s) \leq
-CM(s)^{\frac{d+1}{d}}, \end{equation}
 The argument is valid for
all sufficiently large $M(s),$ which is the only situation we need
to consider provided $A$ was chosen large enough. The same bound
holds for any point $x_0$ where $M(s)$ is attained and by continuity
in some neighborhoods of such points. So, we have \eqref{infin1} in
some open set $U$. Due to smoothness of $\psi$, away from $U$ we
have
$$
\max\limits_{x\not\in U}|\psi(x,\tau)|< M(\tau)
$$
for every $\tau$ during some period of time after $s.$  Thus we
obtain that
\begin{equation}\label{infin}
\frac{d}{ds}M(s) \leq -C M^{\frac{d+1}{d}}(s), \,\,\,M(0) \leq
Ar^{-d}.
\end{equation}
This is valid for all times while $M(s)$ remains sufficiently
large. Solving \eqref{infin}, we get an estimate
\[
M(s) \leq \frac{M(0)}{(1+CM(0)^{1/d}s)^d} \leq
Ar^{-d}(1-CA^{1/d}r^{-1}s)
\]
for all sufficiently small $s.$ This implies
\begin{equation}\label{infinlast}
 \|\psi(\cdot, s)\|_{L^\infty} \leq Ar^{-d}(1-CA^{1/d}r^{-1}s),
\end{equation}
for all sufficiently small $s \leq \gamma(A,d)r.$ Observe that
$\gamma$ is independent of $\psi$ or $v$ other than through the
value of $A,$ which will be chosen below depending on the value of
$B$ only. The estimate \eqref{infinlast} agrees with the properties
of the $(1-\frac{\delta s}{r})\cU_{r+Ks}(\Tm^d)$ class provided that
\begin{equation}\label{linftycon11}
\delta + dK \leq CA^{1/d}.
\end{equation}

 Next, we consider {\bf the concentration condition }
$\int_{\Tm^d}\Omega(x-x_0)|\psi(x)|\,dx \leq r^{1/2}.$ Consider
$x(s) \in \Tm^d$ satisfying
\begin{equation}\label{concdyn}
x'(s) = \overline{v}_{B_r(x(s))} \equiv
\frac{1}{|B_r|}\int_{B_r(x(s))} v(y,s)\,dy, \,\,\,x(0)=x_0.
\end{equation}
Here $B_r(x)$ stands for the ball of radius $r$ centered at $x,$ and
$|B_r|$ is its volume. We will estimate
$\int_{\Tm^d}\Omega(x-x(s))|\psi(x,s)|\,dx.$
 Let us write $\psi(x)=\psi_+(x)-\psi_-(x),$
where $\psi_\pm(x) \geq 0$ and have disjoint support. Let us denote
$\psi_\pm(x,s)$ the solutions of \eqref{psieq} with
$\psi_\pm(x,0)=\psi_\pm(x).$ Then due to linearity and maximum
principle, $|\psi(x,s)| = |\psi_+(x,s)-\psi_-(x,s)| \leq
\psi_+(x,s)+\psi_-(x,s),$ and so
\begin{equation}\label{psipm}
\int_{\Tm^d} \Omega(x-x(s))|\psi(x,s)|\,dx \leq \int_{\Tm^d}
\Omega(x-x(s))\psi_+(x,s)\,dx + \int_{\Tm^d}
\Omega(x-x(s))\psi_-(x,s)\,dx.
\end{equation}
Let us estimate the first integral on the right hand side of
\eqref{psipm}, the second can be handled the same way. We have
\begin{eqnarray}\nonumber
\left|\partial_s \int_{\Tm^d} \Omega(x-x(s)) \psi_+\,dx \right|=
\\ \nonumber \left|\int_{\Tm^d} \left( \Omega(x-x(s))
\left((-v \cdot \nabla)\psi_+ - (-\Delta)^{1/2} \psi_+ \right)-
\nabla(\Omega(x-x(s)))\cdot x'(s) \psi_+ \right)\,dx  \right|= \\
\nonumber \left|\int_{\Tm^d} \nabla(\Omega(x-x(s))) \cdot \left( v -
\overline{v}_{B_r(x(s))} \right)\psi_+\,dx -
\int_{\Tm^d} (-\Delta)^{1/2} \Omega(x-x(s)) \psi_+\,dx \right| \\
\leq C\left(\int_{\Tm^d} |x-x(s)|^{-1/2}|v-
\overline{v}_{B_r(x(s))}||\psi_+|\,dx + \int_{\Tm^d}
|x-x(s)|^{-1/2}|\psi_+|\,dx \right).\label{psiderest}
\end{eqnarray}
We used the divergence free condition on $v$ and \eqref{concdyn} in
the second step, and estimated $|\nabla\Omega(x-x_0)| \leq
C|x-x_0|^{-1/2},$ $|(-\Delta)^{1/2}\Omega(x-x_0)|\leq
C|x-x_0|^{1/2}.$ Let us consider the two integrals in
\eqref{psiderest}. Since $\|\psi_+\|_{L^1}\leq 1/2$ and
$\|\psi_+\|_{L^\infty} \leq Ar^{-d},$ the integral $\int_{\Tm^d}
|x-x(s)|^{-1/2}|\psi_+|\,dx$ is maximal when $\psi_+$ is a
characteristic function of a ball centered at $x(s)$ of radius
$crA^{-1/d}.$ This gives an upper bound of $Cr^{-1/2}A^{1/2d}$ for
this integral. To estimate the first
integral in \eqref{psiderest}, 
split $\Tm^d = \cup_{k=0}^N E_k,$ where
\[ E_k = \{ x:\,\,r2^{k-1} < |x-x(s)| \leq r2^k  \} \cap
\Tm^d,\,\,\,k>0, \,\,\,E_0=B_r(x(s)). \] Recall (see e.g.
\cite{Stein}) that for any BMO function $f$, any ball $B,$ and any
$1 \leq p < \infty,$
\begin{equation}\label{BMOest123}
\|f-\overline{f}_B\|_{L^p(B)} \leq c_p |B|^{1/p}\|f\|_{BMO}.
\end{equation}
 By H\"older's inequality,
\begin{eqnarray*} \int_{B_r(x(s))}
|x-x(s)|^{-1/2}|v-\overline{v}_{B_r(x(s))}||\psi_+|\,dx \leq \\
\||x-x(s)|^{-1/2}\|_{L^{p}(B_r(x(s))}\|v-\overline{v}_{B_r(x(s))}\|_{L^{z}(B_r(x(s)))}\|\psi_+\|_{L^{q}(B_r(x(s)))},
\end{eqnarray*}
where $p^{-1}+z^{-1}+q^{-1}=1.$ Now
\[ \|\psi_+\|_{L^q} \leq
\|\psi_+\|_{L^1}^{1/q}\|\psi_+\|_{L^\infty}^{1-1/q} \leq
A^{1-\frac{1}{q}}r^{\frac{d}{q}-d}.\] Using \eqref{BMOest123}, we
also see that
\[ \|v-\overline{v}_{B_r(x(s))}\|_{L^{z}(B_r(x(s)))} \leq
C(z,d)r^{\frac{d}{z}}B. \] Finally, for any $p < 2d,$
\[ \||x-x(s)|^{-1/2}\|_{L^{p}(B_r(x(s))} \leq
C(p,d)r^{\frac{d}{p}-\frac{1}{2}}. \] Taking $z$ very large, and $p$
very close to $2d,$ we find that for any $q > \frac{2d}{2d-1},$ we
have
\begin{equation}\label{bmoest1}
\int_{B_r(x(s))}
|x-x(s)|^{-1/2}|v-\overline{v}_{B_r(x(s))}||\psi_+(x)|\,dx \leq
CBA^{1-\frac{1}{q}}r^{\frac{d}{q}+\frac{d}{z}+\frac{d}{p}-d-1/2}
\leq C(\sigma,d) B A^\sigma r^{-1/2},
\end{equation}
where $\sigma$ is any number greater than $\frac{1}{2d}.$
Furthermore, for $k>0,$
\begin{eqnarray}
\int_{E_k} |x-x(s)|^{-1/2}|v-\overline{v}_{B_r(x(s))}||\psi_+|\,dx
\leq
C2^{-k/2}r^{-1/2}\int_{E_k} |v-\overline{v}_{B_r(x(s))}||\psi_+(x)|\,dx \leq \nonumber \\
C2^{-k/2}r^{-1/2}\left( \int_{B_{r2^k }(x(s))}
|v-\overline{v}_{B_{r2^{k}}(x(s))}||\psi_+(x)|\,dx+\int_{B_{r2^k
}(x(s))}
|\overline{v}_{B_{r2^{k}}(x(s))}-\overline{v}_{B_r(x(s))}||\psi_+(x)|\,dx
\right). \label{bmoest2}
\end{eqnarray}
Recall that (see, e.g., \cite{Stein})
\[ |\overline{v}_{B_{r2^{k}}(x(s))}-\overline{v}_{B_r(x(s))}|
\leq Ck\|v\|_{BMO}. \] Therefore the last integral in
\eqref{bmoest2} does not exceed $CkB.$ The first integral can be
estimated by
\[ \|v-\overline{v}_{B_{r2^k }}\|_{L^{\frac{q}{q-1}}(B_{r2^k })}\|\psi_+\|_{L^{q}(B_{r2^{k} })}
\leq C(q,d)B2^{k(d-\frac{d}{q})}A^{1-\frac{1}{q}}, \] where $q$ is
any number greater than $1.$ Thus in particular
\begin{equation}\label{bmoest3}
\int_{E_k}
|x-x(s)|^{-1/2}|v-\overline{v}_{B_1(x(s))}||\psi_+(x)|\,dx \leq
CB2^{-3k/8}B(k2^{-k/8}+A^{1/8d})r^{-1/2}
\end{equation}
if $q=\frac{8d}{8d-1}.$ Adding \eqref{bmoest1} and \eqref{bmoest3},
we obtain
\[ \int_{\Tm^d} |x-x(s)|^{-1/2}|v-\overline{v}_{B_r(x(s))}||\psi_+(x,s)|\,dx \leq
CBA^{3/4d}r^{-1/2},\] provided that $A$ is large enough (the
exponent for $A$ can be anything greater then $\frac{1}{2d}).$
Coming back to \eqref{psiderest} and \eqref{psipm}, we see that
\begin{equation}\label{concest11a} \int_{\Tm^d} |x-x(s)|^{1/2}|\psi(x,s)|\,dx
 \leq r^{1/2}+Csr^{-1/2}(A^{1/2d}+BA^{3/4d}). \end{equation}
 This is consistent
with the $(1-\frac{\delta s}{r})\cU_{r+Ks}(\Tm^d)$ class if
\[ \left(1-\frac{\delta s}{r} \right) \left(r+Ks\right)^{1/2} \geq
r^{1/2}+Csr^{-1/2}(A^{1/2d}+BA^{3/4d}). \] Provided that $\gamma$ is
chosen sufficiently small, this condition reduces to
\begin{equation}\label{keyest676}
\frac{K}{2} - \delta > C(A^{1/2d}+BA^{3/4d}).
\end{equation}

Finally, we consider {\bf the $L^1$ norm}. Recall (see e.g.
\cite{CC}) that for a $C^\infty$ function $\psi(x),$
\begin{equation}\label{fracdiff}
(-\Delta)^{1/2}\psi(x) = \lim_{\epsilon \rightarrow 0}\sum\limits_{n
\in \Zm^d} \int\limits_{\Tm^d \cap
|x-y|\geq \epsilon}\frac{\psi(x)-\psi(y)}{|x-y-n|^{d+1}}\,dy.
\end{equation}
Let $S$ be the set where $\psi(x,s)=0,$ and define
\[ D_\pm = \{ x\in \Tm^d \left| \,\pm\psi(x,s)>0 \}. \right. \]
The sets $S$ and $D_\pm$ depend on $s,$ but we will omit this in
notation to save space. Due to \eqref{psieq} and incompressibility
of $v,$ we have
\begin{eqnarray}\nonumber
\partial_s \|\psi(\cdot,s)\|_{L^1} = \int_{\Tm^d \setminus S}
\frac{\psi(x,s)}{|\psi(x,s)|}\left( -v\cdot\nabla\psi(x,s) -
(-\Delta)^{1/2}\psi(x,s)\right)\,dx + \\
\int_{S}|(-\Delta)^{1/2}\psi(x,s)|\,dx  = \label{L1normeq}
-\int_{\Tm^d \setminus S}
\frac{\psi(x,s)}{|\psi(x,s)|}(-\Delta)^{1/2}\psi(x,s)\,dx +
\int_{S}|(-\Delta)^{1/2}\psi(x,s)|\,dx.
\end{eqnarray}
The integral over $S$ is of course nonzero only if the Lebesgue
measure
of $S$ is positive. 
Substituting \eqref{fracdiff} into \eqref{L1normeq} and symmetrizing
with respect to $x,y$ we get
\begin{eqnarray}\nonumber
\partial_s \|\psi(\cdot,s)\|_{L^1} = \\ - \frac12\lim_{\epsilon
\rightarrow 0} \int\limits_{((\Tm^d \setminus S) \times (\Tm^d
\setminus S)) \cap |x-y| \geq \epsilon} \left(
\frac{\psi(x,s)}{|\psi(x,s)|} - \frac{\psi(y,s)}{|\psi(y,s)|}\right)
\sum\limits_{n \in \Zm^d}
\frac{\psi(x,s)-\psi(y,s)}{|x-y-n|^{d+1}}\,dx dy \label{L1nb} \\
- \int_{\Tm^d \setminus S} \frac{\psi(x,s)}{|\psi(x,s)|}
\left(\sum\limits_{n \in \Zm^d} \int_S
\frac{\psi(x,s)}{|x-y-n|^{d+1}}\,dy \right)dx + \int_S \left|
\sum\limits_{n \in \Zm^d} \int_{\Tm^d \setminus S}
\frac{\psi(y,s)}{|x-y-n|^{d+1}}\,dy \right| \,dx. \nonumber
\end{eqnarray}
Observe that the expression under the first integral in \eqref{L1nb}
is non-negative for all $x,y,$ and it is positive if $\psi(x,s)$ and
$\psi(y,s)$ have different signs. Also, observe that
\[ \int_S \left|
\sum\limits_{n \in \Zm^d} \int_{\Tm^d \setminus S}
\frac{\psi(y,s)}{|x-y-n|^{d+1}}\,dy \right| \,dx \leq \int_S
\sum\limits_{n \in \Zm^d} \left|
\int_{D_+}\frac{\psi(y,s)}{|x-y-n|^{d+1}}\,dy +
\int_{D_-}\frac{\psi(y,s)}{|x-y-n|^{d+1}}\,dy \right|. \] Therefore,
the combined contribution of the last line in \eqref{L1nb} over
every cell is less than or equal to zero. Leaving only the central
cell contributions in \eqref{L1nb}, we get
\begin{eqnarray}
\nonumber\partial_s \|\psi(\cdot,s)\|_{L^1} \leq  \\
\label{L1bnd11} -\int\limits_{D_+} \psi(x,s) \int\limits_{D_-}
\frac{dy}{|x-y|^{d+1}}\,dydx + \int\limits_{D_-} \psi(y,s)
\int\limits_{D_+} \frac{dx}{|x-y|^{d+1}}\,dxdy - \\
\nonumber - \int_{D_+ \cup D_-} |\psi(x,s)| \left( \int_S
\frac{1}{|x-y|^{d+1}}\,dy \right)dx+\int_S \left|
\int_{D_+}\frac{\psi(y,s)}{|x-y|^{d+1}}\,dy +
\int_{D_-}\frac{\psi(y,s)}{|x-y|^{d+1}}\,dy \right|.
\end{eqnarray}
Without loss of generality, we can assume that $1 \geq
\|\psi(\cdot,s)\|_{L^1} \geq 9/10$ for every $s$ we consider, since
otherwise the $L^1$ condition is already satisfied. Also, due to
\eqref{concest11a} we can assume that $\int_{\Tm^d}
\Omega(x-x(s))|\psi(x,s)|\,dx \leq \frac{11}{10}r^{1/2}$ provided
that the time interval $[0,\gamma r]$ that we consider is
sufficiently small, with $\gamma = \gamma(A,B).$ These two bounds
imply that $\int_{\Tm^d \cap |x-x(s)|\leq
400r} |\psi(x,s)|\,dx \geq 4/5.$ 
The mean zero condition leads to
\begin{equation}\label{mzconest}
\pm\int_{D_\pm \cap \{|x-x(s)| \leq 400r\}} \psi(x,s)\,dx \geq 3/10.
\end{equation}
Let us denote $\widetilde{D}_\pm = D_\pm \cap \{|x-x(s)| \leq
400r\},$ $\widetilde{S} = S \cap \{|x-x(s)| \leq 400r\}.$ Observe
that if $x \in \widetilde{S},$ then, by \eqref{mzconest},
\[ \pm \int_{D_\pm}\frac{\psi(y,s)}{|x-y|^{d+1}}\,dy \geq
\pm \int_{\widetilde{D}_\pm}\frac{\psi(y,s)}{|x-y|^{d+1}}\,dy \geq C
r^{-d-1}. \] This implies that due to cancelation in the last term
of \eqref{L1bnd11}, we can estimate the last line of \eqref{L1bnd11}
from above by $- C|\widetilde{S}|r^{-d-1}.$ Reducing the integration
in the second line of \eqref{L1bnd11} to $\widetilde{D}_\pm,$ we
obtain
\begin{equation}\label{finL1est}
\partial_s \|\psi(\cdot,s)\|_{L^1} \leq -C r^{-d-1}\left(|\widetilde{D}_-|\int_{\widetilde{D}_+}
\psi(x,s)\,dx  + |\widetilde{D}_+|\int_{\widetilde{D}_-}
\psi(x,s)\,dx +|\widetilde{S}|\right) \leq -cr^{-1},
\end{equation} where $c$ is a fixed positive constant.
Here in the last step we used \eqref{mzconest} and
$|\widetilde{D}_+| + |\widetilde{D}_-| +|\widetilde S| \geq C r^d.$
The estimate \eqref{finL1est} is consistent with $(1-\frac{\delta
s}{r})\cU_{r+Ks}(\Tm^d)$ class if $\delta \leq c.$

It remains to observe that, if $A=A(B,d)$ is chosen sufficiently
large, one can indeed find $K$ and $\delta$ so that the conditions
\eqref{linftycon11}, \eqref{keyest676} and the $\delta \leq c$
condition arising from the $L^1$ norm estimate are all satisfied. It
is also clear from the proof that \eqref{finreslem} then holds for
all $s\leq \gamma(B,d) r.$ The only restriction from above on the
value of $r$ comes from the $L^\infty$ norm condition, which has to
be consistent with $L^1$ and concentration conditions. For
convenience, we chose to cap the value of $r$ at $1.$
\end{proof}

The proof of Theorem~\ref{key} is now straightforward.
\begin{proof}[Proof of Theorem~\ref{key}]
From Lemma~\ref{keylemma} it follows that for any $s>0,$ $\psi(x,s)
\in f(s) \cU_{r+Ks}(\Tm^d)$ provided that $f'(s) \geq -
\frac{\delta}{r+Ks}f(s).$ Solving this differential equation, we
obtain that the factor $f(s) = \left(\frac{r}{r+Ks}
\right)^{\delta/K}$ is acceptable.
\end{proof}

\section{The Critical SQG equation and further discussion}\label{SQG}

Theorem~\ref{main1} provides an alternative path to the proof of
existence of global regular solutions to the critical surface
quasi-geostrophic equation:
\begin{equation}
\left\{ \aligned &\theta_t=u\cdot\nabla\theta
-(-\Delta)^{1/2}\theta, \,\,\,\,\theta(x,0)=\theta_0(x),
\\ \nonumber
& u=(u_1,u_2)=(-R_2\theta, R_1\theta),
\endaligned
\right.
\end{equation}
where $\theta\,:\, \Rm^2\to  \Rm^2$ is a periodic scalar function,
and $R_1$ and $R_2$ are the usual Riesz transforms in $\Rm^2$.
 Indeed, the local existence and uniqueness of smooth solution starting
from $H^1$ periodic initial data is known (see e.g. \cite{dong}).
The $L^\infty$ norm of the solution does not increase due to the
maximum principle (see e.g. \cite{CC}), which implies uniform bound
on the BMO norm of the velocity. Since the local solution is smooth,
one can apply Theorem~\ref{main1}. This, similarly to \cite{CV},
implies a uniform bound on some H\"older norm of the solution
$\theta.$ This improvement over the $L^\infty$ control is sufficient
to show the global regularity (see \cite{CV} or \cite{CW} for
slightly different settings which can be adapted to our case in a
standard way).

One can pursue a number of generalizations of Theorem~\ref{main1},
for instance reducing assumptions on smoothness of solution,
velocity, or initial data. However we chose to present here the case
with the most transparent proof containing the heart of the matter.
As follows from the proof, the role of the BMO space is mainly the
right scaling: the BMO is the most general function space for which
\eqref{BMOest123} is available. The BMO scaling properties are of
course also crucial for the proof of \cite{CV} to work.

\noindent {\bf Acknowledgement.} \rm Research of AK has been
supported in part by the NSF-DMS grant 0653813. Research of FN has
been partially supported by the NSF-DMS grant 0501067. AK thanks for
hospitality the Department of Mathematics of the University of
Chicago, where part of this work was carried out.

\end{document}